\documentclass[oneside,english]{amsart}
\usepackage[T1]{fontenc}
\usepackage[latin1]{inputenc}
\usepackage{amssymb}
\usepackage{graphicx}

\makeatletter
  \theoremstyle{plain}
  \newtheorem{prop}{Proposition}
 \theoremstyle{definition}
 \newtheorem{defn}{Definition}
  \theoremstyle{remark}
  \newtheorem*{rem*}{Remark}
  \theoremstyle{plain}
  \newtheorem{lem}{Lemma}
  \theoremstyle{plain}
  \newtheorem*{conjecture*}{Conjecture}
  \theoremstyle{definition}
  \newtheorem*{example*}{Example}

\usepackage{a4wide}

\usepackage{babel}
\makeatother
\begin{document}
\newcommand{\cat}[1]{\mathsf{#1}}

\newcommand{\map}{\rightarrow}

\newcommand{\Hom}{\mathrm{Hom}}

\newcommand{\HOM}{\mathbf{Hom}}

\newcommand{\rodd}[1]{\mathbb{R}^{0|#1}}

\title{Differential worms and generalized manifolds}

\author{Pavol \v Severa}

\curraddr{Section de Math\'ematiques, Universit\'e de Gen\`eve, rue du Li\`evre
2-4, 1211 Geneva, Switzerland, on leave from Dept.~of Theoretical
Physics, Mlynsk\'a dolina F2, 84248 Bratislava, Slovakia}

\email{pavol.severa@gmail.com}

\thanks{partially supported by the Swiss National Science Foundation}

\maketitle

\begin{abstract}
We study differential forms and their higher-order generalizations
by interpreting them as functions on map spaces. We get a series of
approximations of {}``generalized manifolds'' (i.e.~of sheaves
and stacks) somewhat akin to Taylor series.
\end{abstract}
\tableofcontents{}

\section{Multidifferential algebras on manifolds}

We shall make some computations with natural multidifferential algebras
that generalize the algebra of differential forms on a manifold. The
elements of these algebras will be called differential worms. In this
section we just want to see how to compute in local coordinates and
what is the behaviour under coordinate changes and pullbacks. 

Let us start with an open subset $U$ of $\mathbb{R}^{n}$ (since
we are going to compute in local coordinates). For any $k\in\mathbb{N}$
we shall construct a $k$-fold differential graded-commutative algebra
$\Omega_{[k]}(U)$ (i.e.~a $k$-fold complex with a compatible algebra
structure), generalizing differential forms on $U$; $\Omega_{[1]}(U)$
will be just $\Omega(U)$. 

The algebra $\Omega_{[k]}(U)$ of \emph{differential worms of level}
$k$ \emph{}will have the following properties: its zero-degree subalgebra
will be $C^{\infty}(U)$ and for every function $f\in C^{\infty}(U)$
and every differential $d_{a}$ among the $k$ differentials $d_{1},\dots,d_{k}$
on $\Omega_{[k]}(U)$, we will have\begin{equation}
d_{a}f=\frac{\partial f}{\partial x^{i}}\, d_{a}x^{i}.\label{eq:dif}\end{equation}
Out of it we easily compute\begin{equation}
d_{a}d_{b}f=\frac{\partial f}{\partial x^{i}}\ d_{a}d_{b}x^{i}+\frac{\partial^{2}f}{\partial x^{i}\partial x^{j}}\ d_{a}x^{i}\ d_{b}x^{j},\label{eq:diff}\end{equation}
\begin{eqnarray}
d_{a}d_{b}d_{c}f & = & \frac{\partial f}{\partial x^{i}}\ d_{a}d_{b}d_{c}x^{i}+\frac{\partial^{2}f}{\partial x^{i}\partial x^{j}}\left(d_{a}x^{i}\ d_{b}d_{c}x^{j}+d_{b}x^{i}\ d_{c}d_{a}x^{j}+d_{c}x^{i}\ d_{a}d_{b}x^{j}\right)+\nonumber \\
 &  & {}+\frac{\partial^{3}f}{\partial x^{i}\partial x^{j}\partial x^{k}}\ d_{a}x^{i}d_{b}x^{j}d_{c}x^{k}\label{eq:difff}\end{eqnarray}
and so on.

The algebra $\Omega_{[k]}(U)$ is defined as the graded-commutative
algebra freely generated by the algebra $C^{\infty}(U)$ and by elements
$d_{a_{1}}d_{a_{2}}\dots d_{a_{p}}x^{i}$, $1\leq a_{1}<a_{2}<\dots<a_{p}\leq k$,
$1\leq i\leq n$.

We defined $\Omega_{[k]}(U)$ using coordinates, but it is easy to
see what happens when we pass to a different system of coordinates
$\tilde{x}^{i}$. First of all, the relation (\ref{eq:dif}) is valid
in any coordinates, since \[
\frac{\partial f}{\partial\tilde{x}^{j}}\, d_{a}\tilde{x}^{j}=\frac{\partial f}{\partial\tilde{x}^{j}}\frac{\partial\tilde{x}^{j}}{\partial x^{i}}\, d_{a}x^{i}=\frac{\partial f}{\partial x^{i}}\, d_{a}x^{i}=d_{a}f.\]
 Hence also the relations (\ref{eq:diff}), (\ref{eq:difff}) etc.~hold
in any coordinates. Now we can pass between any two systems of coordinates
-- we just set the $f$ on the LHS of (1), (2), (3) etc.~to be coordinates
from one system and the $x$'s on the RHS to be coordinates from the
other system. Since we can pass between systems of coordinates, we
can also define $\Omega_{[k]}(M)$ on any manifold (it is easy to
see that this definition is consistent). Moreover, for any smooth
map $M\map N$ we have a pullback $\Omega_{[k]}(N)\map\Omega_{[k]}(M)$,
again given in coordinates by (1), (2), (3) etc. We should notice
that the algebra $\Omega_{[k]}(M)$ belongs to higher-orded geometry
if $k\geq2$, since the transition formulas contain up to $k$'th
derivatives.

The definition we just presented was very low-brow. One could give
a more natural (coordinate-free) definition, using a universal property
satisfied by $\Omega_{[k]}(M)$. We shall, however, wait until Section
\ref{sec:forms-worms}, where $\Omega_{[k]}(M)$'s will get a geometrical
explanation, together with a much richer structure than that of a
$k$-fold differential algebra.

\section{Late informal introduction\label{sec:Late-introduction}}

It is well known that one can see differential forms as functions
on a supermanifold and de Rham differential as a vector field on the
supermanifold. It is however less well known that the supermanifold
and the vector field have a very simple geometrical explanation. Namely,
differential forms on $M$ are functions on the supermanifold $\HOM(\rodd{1},M)$
of all maps $\rodd{1}\map M$ (here $\rodd{1}$ denotes the odd line),
while both the differential and the degrees of differential forms,
i.e.~the structure of a complex on $\Omega(M)$, come from the action
of $\mathit{Diff}(\rodd{1})$ on $\HOM(\rodd{1},M)$. We found this
geometrical explanation in \cite{kon}.

It is then straightforward to see $\Omega_{[k]}(M)$ as functions
on $\HOM(\rodd{k},M)$. This gives us an action of $\mathit{Diff}(\rodd{k})$
on $\Omega_{[k]}(M)$, which is (for $k\geq2$) a stronger structure
than just a $k$-fold complex. We thus get not just a geometric meaning
of differential worms, but also some new interesting properties. This
point of view of worms also suggests the next developpement, as follows.

Besides differential forms there are other interesting differential
graded-commutative algebras (DGCAs). For example, if $\mathfrak{g}$
is a Lie algebra, on the Grassmann algebra $\bigwedge\mathfrak{g}^{*}$
there is the Chevalley-Eileinberg differential that makes it to a
DGCA. As we'll see, these other DGCAs can also be seen as differential
forms on some {}``generalized manifolds''. By differential forms
we mean functions on the space of maps from $\rodd{1}$ to the generalized
manifold. The generalized manifolds we have in mind are contravariant
functors from the category of (super)manifolds. If $F$ is such a
functor, the idea is to see $F(M)$ as the {}``space of all maps
$M\map F$'', and thus differential forms on $F$ are functions on
$F(\rodd{1})$. This way we can construct $\Omega_{[k]}(F)$ for any
$k$, and moreover, approximate $F$ using $\Omega_{[k]}(F)$ and
the action of $\mathit{Diff}(\rodd{k})$ on $\Omega_{[k]}(M)$ (the
approximation is $\mathrm{app}_{k}F(M)=\{\mathit{Diff}(\rodd{k})\mbox{-equivariant algebra morphisms }\Omega_{[k]}(F)\map\Omega_{[k]}(M)\}$).
We thus get a sequence of approximations of $F$, akin to Taylor series.

As a simple example, let $G$ be a Lie group, and let us define $F(M)=\HOM(M,G)/G$.
In this case $\Omega(F)$ is the Chevalley-Eilenberg DGCA $\bigwedge\mathfrak{g}^{*}$.
The approximation $\mathrm{app}_{1}F$ is given by \[
\mathrm{app}_{1}F(M)=\{\mbox{DGCA maps }{\textstyle \bigwedge}\mathfrak{g}^{*}\map\Omega(M)\}=\{\mbox{flat connections on }M\}.\]
All higher approximations of $F$ are in this case equal to $\mathrm{app}_{1}F$.

This paper is inspired, besides the fact that differential forms are
functions on $\HOM(\rodd{1},M)$, by Sullivan's Rational homotopy
theory \cite{sull}, where one interprets general (reasonable) DGCAs
as if they were differential forms on manifolds. In particular, one
considers DGCA maps $A\map\Omega(M)$ as maps from $M$ to some generalized
space (out from this one can e.g.~construct a simplicial set by taking
$M$'s to be simplices of all dimensions, and then get the homotopy
type of the generalized space as the geometric realization). We are
interested in the inverse operation --- getting a DGCA $\Omega(F)$
out of a generalized space $F$, and also getting the algebras $\Omega_{[k]}(F)$
that approximate $F$ better than $\Omega(F)$ does.

Most of the ideas of this paper appeared in the unpublished work \cite{ks};
here we attempt to state them in a formally correct way. Differential
worms were introduced independently (and with a saner name) by Vinogradov
and Vitagliano \cite{vv}, from a different point of view.

\section{Differential forms and differential worms as functions on map spaces\label{sec:forms-worms}}

The basic hero of this paper is the supermanifold $\HOM(\rodd{k},X)$
of all maps $\rodd{k}\map X$, where $X$ is any (super)manifold.
It is defined, up to unique isomorphisms, by the following property:
for any supermanifold $Y$, a map $Y\map\HOM(\rodd{k},X)$ is the
same as a map $Y\times\rodd{k}\map X$.

It is easy to see that $\HOM(\rodd{k},X)$ exists for any $k$ and
$X$, and moreover, is naturally (i.e.~functorially in $X$) isomorphic
with something well-known:

\begin{prop}
\label{pro:basic}For any supermanifold $X$, $\HOM(\rodd{1},X)$
is naturally isomorphic to $\Pi TX$, and thus (using induction) $\HOM(\rodd{k},X)$
is naturally isomorphic to $(\Pi T)^{k}X$.
\end{prop}
Here $\Pi TX$ is the odd tangent bundle of $X$. Recall that differential
forms on $X$ are functions on $\Pi TX$.%
\footnote{more precisely, functions that are polynomial on the fibres of $\Pi TX$;
general functions on $\Pi TX$ are called pseudodifferential forms
on $X$. If $X$ is a manifold then every pseudodifferential form
is actually a differential form.%
} The proof of the proposition is straightforward: let $\theta$ be
the coordinate on $\rodd{1}$ and $x^{i}$ be local coordinates on
$X$, so that $x^{i}$, $dx^{i}$ are local coordinates on $\Pi TX$.
If $\rodd{1}\map X$ is a map parametrized by some $Y$ (that is,
a map $\rodd{1}\times Y\map X$), expanding to Taylor series in $\theta$
it is of the form $x^{i}(\theta)=x^{i}(0)+\xi^{i}\theta$; identifying
$x^{i}(0)$ with $x^{i}$ and $\xi^{i}$ with $dx^{i}$ we get the
(local) isomorphism between $\HOM(\rodd{1},X)$ and $\Pi TX$. This
identification is independent of the choice of coordinates $x^{i}$,
since using other coordinates $\tilde{x}^{i}$ we have \[
\tilde{x}(x(\theta))=\tilde{x}(x(0)+\xi\theta)=\tilde{x}(x(0))+\frac{\partial\tilde{x}}{\partial x}\xi\theta,\]
i.e.~$\tilde{\xi}=\frac{\partial\tilde{x}}{\partial x}\xi$, just
as $d\tilde{x}=\frac{\partial\tilde{x}}{\partial x}dx$.

The isomorphism of $\HOM(\rodd{1},X)$ with $\Pi TX$ is certainly
not surprising: the relation $\theta^{2}=0$ basically says that maps
$\rodd{1}\map X$ are 1-jets of curves in $X$, i.e.~tangent vectors
in $X$. It gives us, however, an interesting explanation of the deRham
differential (which, being a derivation on $\Omega(X)$, is a vector
field on $\Pi TX$). Namely, on $\HOM(\rodd{1},X)$ we have a right
action of the supersemigroup $\HOM(\rodd{1},\rodd{1})$, and this
action gives us the structure of a complex on $\Omega(X)$ (a left
action of $\HOM(\rodd{1},\rodd{1})$ on a vector space $V$ is the
same as a structure of a non-negatively graded complex on $V$).

Let us compute the action of $\HOM(\rodd{1},\rodd{1})$ explicitly.
Given a transformation $\theta\mapsto\theta'=a\theta+\beta$ of $\rodd{1}$,
we have\[
x^{i}(\theta')=x^{i}(a\theta+\beta)=\left(x^{i}(0)+\xi^{i}\beta\right)+a\xi^{i}\theta,\]
i.e.\[
x^{i}\mapsto x^{i}+dx^{i}\,\beta,\qquad dx^{i}\mapsto a\, dx^{i}.\]
If we take infinitesimal generators of $\HOM(\rodd{1},\rodd{1})$,
$\partial/\partial\theta$ and $\theta\partial/\partial\theta$, we
get the following:

\begin{prop}
The vector fields $\partial/\partial\theta$ and $\theta\,\partial/\partial\theta$
on $\rodd{1}$ act on $\Pi TM$ as\[
dx^{i}\frac{\partial}{\partial x^{i}}\quad\textrm{and}\quad dx^{i}\frac{\partial}{\partial(dx^{i})},\]
i.e.~as the de Rham differential and the degree.
\end{prop}
Now we can make an obvious generalization, by passing from $\rodd{1}$
to $\rodd{k}$. On $\HOM(\rodd{k},X)$ we have a right action of $\HOM(\rodd{k},\rodd{k})$.
In particular, the vector fields $\partial/\partial\theta^{a}$ on
$\rodd{k}$ (where $\theta^{a}$s, $1\leq a\leq k$, are coordinates
on $\rodd{k}$) give rise to $k$ odd vector fields $d_{a}$ on $\HOM(\rodd{k},X)$,
satisfying $d_{a}d_{b}+d_{b}d_{a}=0$, i.e.~to $k$ anticommuting
differentials. In local coordinates $x^{i}$ on $X$, a map $\rodd{k}\map X$
is given by functions $x^{i}(\theta)$ on $\rodd{k}$ (we suppress
the inessential dependence on $Y$), and we have\[
x^{i}(\theta)=\exp(\theta^{a}d_{a})\, x^{i}(0).\]
We denote $x^{i}(0)$ simply by $x^{i}$ (a slight abuse of notation:
$x^{i}$'s were functions on $X$ and now they become functions on
$\HOM(\rodd{k},X)$). As an example, for $k=2$ we get \[
x^{i}(\theta^{1},\theta^{2})=x^{i}+\theta^{1}d_{1}x^{i}+\theta^{2}d_{2}x^{i}+\theta^{2}\theta^{1}d_{1}d_{2}x^{i}.\]

Now we see that \[
d_{a_{1}}d_{a_{2}}\dots d_{a_{p}}x^{i}\quad(1\leq a_{1}<a_{2}<\dots<a_{p}\leq k,\,0\leq p\leq k)\]
are coordinates on $\HOM(\rodd{k},X)$ and we can identify $\Omega_{[k]}(X)$
with functions on $\HOM(\rodd{k},X)$.%
\footnote{more precisely, $\Omega_{[k]}(X)$ consists of the functions that
are polynomial in the differentials of $x^{i}$'s; more invariantly,
they are the functions that are in finite-dimensional $\HOM(\rodd{k},\rodd{k})$-invariant
subspaces of $C^{\infty}(\HOM(\rodd{k},X))$%
}

$\HOM(\rodd{k},\rodd{k})$ acts on $\HOM(\rodd{k},X)$ from the right,
and therefore on $\Omega_{[k]}(X)$ from the left. This gives $\Omega_{[k]}(X)$
a stronger structure than that of a $k$-differential algebra. For
example when $k=2$, vector fields on $\rodd{2}$ have basis\[
\frac{\partial}{\partial\theta^{a}},\quad\theta^{b}\frac{\partial}{\partial\theta^{a}},\quad\theta^{1}\theta^{2}\frac{\partial}{\partial\theta^{a}},\]
which generates on $\HOM(\rodd{2},X)$ the vector fields\begin{eqnarray*}
d_{a} & = & d_{a}x^{i}\frac{\partial}{\partial x^{i}}+\epsilon_{ab}\, d_{1}d_{2}x^{i}\frac{\partial}{\partial(d_{b}x^{i})}\\
E_{a}^{b} & = & d_{a}x^{i}\frac{\partial}{\partial(d_{b}x^{i})}+\delta_{a}^{b}\, d_{1}d_{2}x^{i}\frac{\partial}{\partial(d_{1}d_{2}x^{i})}\\
R_{a} & = & \epsilon_{ab}\, d_{b}x^{i}\frac{\partial}{\partial(d_{1}d_{2}x^{i})},\end{eqnarray*}
out of which only $d_{1}$, $d_{2}$ (differentials) and $E_{1}^{1}$,
$E_{2}^{2}$ (degrees) take part in the bicomplex structure.

For further properties of worms see Appendix, and also the last example
in Section \ref{sec:Examples} that contains a decomposition of $\Omega_{[2]}(X)$
to indecomposable representations of $\HOM(\rodd{2},\rodd{2})$. Some
other properties can be found in \cite{ks}.

\section{Reminder on presheaves (generalized objects)}

The basic reference for this section is \cite{sga4}, exposé I. Let
$\mathsf{C}$ be a category. A \emph{presheaf} on $\cat{C}$ is a
functor $\cat{C}^{o}\map\cat{Set}$; the category of presheaves%
\footnote{we make the usual hyper-correct assumption of working in some universe
of sets to avoid set-of-all-sets-like problems%
} on $\cat{C}$ (with natural transformations as morphisms) is denoted
$\hat{\cat{C}}$. Presheaves can be reasonably viewed as generalized
objects of $\cat{C}$, with $F(X)$ ($F\in\hat{\cat{C}}$, $X\in\cat{C}$)
interpreted as the set of morphisms $X\map F$. Namely, any object
$Y\in\cat{C}$ gives us a presheaf $Y\in\hat{\cat{C}}$ via $Y(X)=\Hom_{\cat{C}}(X,Y)$.
For any $X\in\cat{C}$ and $F\in\hat{\cat{C}}$ then $\Hom_{\hat{\cat{C}}}(X,F)\cong F(X)$.
This way $\cat{C}$ is identified with a full subcategory of $\hat{\cat{C}}$
(which is the excuse for denoting $X\in\cat{C}$ and the corresponding
presheaf $X\in\hat{\cat{C}}$ by the same letter). 

A presheaf isomorphic to some $X\in\cat{C}$ is said to be \emph{representable}.
For example, if $U$, $V$ are objects of $\cat{C}$, the presheaf
$U\times V$ is defined as $U\times V(X)=U(X)\times V(X)$; if it
is representable, the corresponding object of $\cat{C}$ (defined
up to a unique isomorphism) is called the cartesian product of $U$
and $V$, and denoted (somewhat abusively) $U\times V$ as well. Similarly,
the presheaf $\HOM(U,V)$ is defined by $\HOM(U,V)(X)=\Hom(U\times X,V)$.
If it is representable, the corresponding object of $\cat{C}$ is
called the internal Hom from $U$ to $V$, and is still denoted $\HOM(U,V)$. 

Proposition \ref{pro:basic} says that in the category of supermanifolds,
$\HOM(\rodd{k},X)$ is representable for any $k$ and $X$. It will
be quite conventient to speak of $\HOM(Y,X)$ even when it is not
representable; we just have to keep in mind that it is no longer a
supermanifold, just a presheaf.

\section{Differential forms and presheaves on families of supermanifolds}

Let us denote by $\cat{SM}$ the category of finite-dimensional $C^{\infty}$-supermanifolds. 

Recall from Section \ref{sec:Late-introduction} that we want to define
differential forms on a generalized (super)manifold $F$ as the algebra
of functions on $F(\rodd{1})$. For this to make sense, $F(\rodd{1})$
should be itself a supermanifold, and thus $F$ should be a functor
$\cat{SM}^{o}\map\cat{SM}$ rather than $\cat{SM}^{o}\map\cat{Set}$.
To get interesting examples we would have to admit some infinite-dimensional
supermanifolds to $\cat{SM}$. We prefer not to do it, and rather
use $\widehat{\cat{SM}}$ as a formal replacement of infinite-dimensional
supermanifolds, i.e.~to consider functors $\cat{SM}^{o}\map\widehat{\cat{SM}}$
(which is the same as presheaves on $\cat{SM}\times\cat{SM}$). We
need these contravariant functors to be functorial in this strong
sense: for any two supermanifolds $X$, $Y$ we need a morphism $\HOM(X,Y)\times F(Y)\map F(X)$,
and these morphisms are required to be coherent under the composition
of $\HOM$'s. Such a functor is the same as a presheaf on the following
category, which has the same objects as $\cat{SM}×\cat{SM}$, but
more morphisms.

\begin{defn}
\emph{The category} $\cat{PSM}$ (of product families of supermanifolds)
has pairs of supermanifolds as objects, and morphisms $(X_{1},B_{1})\map(X_{2},B_{2})$
are commutative squares\[
\begin{array}{ccc}
X_{1}×B_{1} & \map & X_{2}×B_{2}\\
\downarrow &  & \downarrow\\
B_{1} & \map & B_{2}\end{array}\]
where the vertical arrows are the canonical projections.
\end{defn}
It means that a morphism $(X_{1},B_{1})\map(X_{2},B_{2})$ is the
same as a map $B_{1}\map B_{2}$ and a map $X_{1}×B_{1}\map X_{2}$,
or in other words, $\Hom\left((X_{1},B_{1}),(X_{2},B_{2})\right)\cong(B_{2}×\HOM(X_{1},X_{2}))(B_{1})$.
This observation gives a morphism $\HOM(X,Y)\times F(Y,\cdot)\map F(X,\cdot)$
in $\widehat{\cat{SM}}$ for any $X,Y\in\cat{SM}$ and $F\in\widehat{\cat{PSM}}$,
in particular, a right action of $\HOM(\rodd{1},\rodd{1})$ on $F(\rodd{1},\cdot)$. 

\begin{defn}
\label{def:forms-on-functor}If $F\in\widehat{\cat{PSM}}$ and $F(\rodd{1},\cdot)\in\widehat{\cat{SM}}$
is representable, \emph{differential forms} on $F$ are functions
on (the supermanifold representing) $F(\rodd{1},\cdot)$. Similarly,
level-$n$ \emph{differential worms} are defined by substituting $\rodd{1}$
with $\rodd{n}$.
\end{defn}
In particular, if $F$ is represented by $(X,\{\mbox{point}\})$ then
$F(\rodd{1},\cdot)$ is represented by the supermanifold $\HOM(\rodd{1},X)=\Pi TX$,
so that differential forms on $F$ are the same as differential forms
on $X$. 

Presheaves on $\cat{PSM}$ are a reasonable definition for the {}``generalized
manifolds'' of Section \ref{sec:Late-introduction}; they will be
the basic object of our study.

\begin{rem*}
Presheaves on the category $\cat{PSM}$ might look a bit exotic. However,
{}``natural'' presheaves on the category of manifolds usually admit
a natural extension to $\cat{PSM}$. Consider e.g.~the presheaf of
differential $k$-forms. Firstly, differential forms make sense on
supermanifolds, so we have a presheaf on $\cat{SM}$. Secondly, if
$X$ and $B$ are supermanifolds, we can consider families of differential
$k$-forms on $X$ smoothly parametrized by $B$; this way we get
a presheaf on $\cat{PSM}$. Roughly speaking, smooth families of something
on $X$, parametrized by $B$, are the presheaves on $\cat{PSM}$
that we are interested in.
\end{rem*}

\section{Reminder on presheaves (continued)}

Here we shall recall the relation between presheaves on a category
and on its full subcategory: if $\cat{C}\subset\cat{D}$ is a full
subcategory then we can see $\hat{\cat{C}}$ as a full subcategory
of $\hat{\cat{D}}$, we have a projection $\mathrm{app}:\hat{\cat{D}}\map\hat{\cat{D}}$
onto $\hat{\cat{C}}$ and a natural transformation $\mathrm{id}_{\hat{\cat{D}}}\map\mathrm{app}$.
It will be the basis for approximations of {}``generalized manifolds''
using their worms (hence the notation {}``app'' for the functor),
as outlined in Section~\ref{sec:Late-introduction}.

Let first $\cat{C}$, $\cat{D}$ be any categories and $u:\cat{C}\map\cat{D}$
any functor. It induces a functor $u^{*}:\hat{\cat{D}}\map\hat{\cat{C}}$
via $u^{*}(F)=F\circ u$. The functor $u^{*}$ admits a right adjoint
$u_{*}:\hat{\cat{C}}\map\hat{\cat{D}}$, that can be defined via $u_{*}(F)(Z)=\Hom_{\hat{\cat{C}}}(u^{*}(Z),F)$
for $F\in\hat{\cat{C}}$ and $Z\in\cat{D}$. If now $\cat{C}$ is
a full subcategory%
\footnote{There is, of course, no real difference between a fully faithful functor
and a full subcategory. Sometimes we shall even commit the crime of
calling a category $\cat{C}$ to be a full subcategory of $\cat{D}$
when all we have is a fully faithful functor $\cat{C}\map\cat{D}$,
provided the functor is clear from the context ($\cat{C}$ is then
just equivalent to a full subcategory of $\cat{D}$).%
} of $\cat{D}$ and $u$ is the inclusion (the only situation we shall
meet; in that case $u^{*}$ is simply the restriction) then $u^{*}\circ u_{*}$
is (isomorphic to) $\mathrm{id}_{\hat{\cat{C}}}$ and consequently
$u_{*}$ is fully faithful. We can thus use $u_{*}$ to identify $\hat{\cat{C}}$
with a full subcategory of $\hat{\cat{D}}$, and $u^{*}$ gives us
a projection $\hat{\cat{D}}\map\hat{\cat{C}}$; the above mentioned
functor $\mathrm{app}:\hat{\cat{D}}\map\hat{\cat{D}}$ will be defined
as $u_{*}\circ u^{*}$ and the morphism $\mathrm{id}_{\hat{\cat{D}}}\map\mathrm{app}$
comes from the fact that $u^{*}$ and $u_{*}$ form an adjoint pair.

We finish with a simple condition that forces $u^{*}$ (and thus $u_{*}$)
to be an equivalence of categories.

\begin{lem}
\label{lem:trivequiv}Suppose that $u:\cat{C}\map\cat{D}$ is a fully
faithful functor, and that any object $Y\in\cat{D}$ is a retract
of some object $X\in\cat{C}$, i.e.~that there are morphisms $Y\map u(X)\map Y$
that compose to $\mathrm{id}_{Y}$. Then $u^{*}$, and consequently
$u_{*}$, is an equivalence of categories.
\end{lem}

\section{Approximations of presheaves\label{sec:Approximations}}

\begin{defn}
For $n\in\mathbb{N}$, \emph{the category} $\cat{PSM}_{n}$ is the
full subcategory of $\cat{PSM}$ with objects $(\rodd{n},B)$, $B\in\cat{SM}$
.
\end{defn}
We have $\Hom\left((\rodd{n},B_{1}),(\rodd{n},B_{2})\right)\simeq\Hom(B_{1},B_{2})×\HOM(\rodd{n},\rodd{n})(B_{1})$,
so we get the following lemma:

\begin{lem}
\label{lem:PSMn}An object of $\widehat{\cat{PSM}_{n}}$ can be equivalently
described as an object of $\widehat{\cat{SM}}$ with a right action
of $\HOM(\rodd{n},\rodd{n})$; a morphism of $\widehat{\cat{PSM}_{n}}$
corresponds to an equivariant morphism of $\widehat{\cat{SM}}$.
\end{lem}
This lemma means that if $F\in\widehat{\cat{PSM}}$ then taking $F(\rodd{n},\cdot)\in\widehat{\cat{SM}}$
together with the right action of $\HOM(\rodd{n},\rodd{n})$ is equivalent
to restricting $F$ from $\cat{PSM}$ to $\cat{PSM}_{n}$. In other
words, if $F(\rodd{n},\cdot)$ happens to be representable, the algebra
of level-$n$ worms on $F$ together with the (left) action of $\HOM(\rodd{n},\rodd{n})$
is equivalent to the restriction of $F$ to $\cat{PSM}_{n}$. In particular,
the differential graded algebra of differential forms on $F$ is equivalent
to the restriction of $F$ to $\cat{PSM}_{1}$.

Next we shall see that restriction to $\cat{PSM}_{k}$ can be restored
from restriction to $\cat{PSM}_{l}$ whenever $k\leq l$. To this
end we define some auxiliary categories:

\begin{defn}
\emph{The category} $\cat{PSM}_{\leq n}$ is the full subcategory
of $\cat{PSM}$ with objects $(X,B)$ such that $X=\rodd{m}$ with
$m\leq n$.
\end{defn}
These categories form a chain $\cat{PSM}_{\leq0}\subset\cat{PSM}_{\leq1}\subset\cat{PSM}_{\leq2}\subset\cdots\subset\cat{PSM}$.
Moreover the inclusion $\cat{PSM}_{n}\subset\cat{PSM}_{\leq n}$ satisfies
the assumption of Lemma \ref{lem:trivequiv}, i.e.

\begin{lem}
Let $u:\cat{PSM}_{n}\map\cat{PSM}_{\leq n}$ be the inclusion. Then
$u^{*}$ is an equivalence of categories.
\end{lem}
This means that there is no difference between restriction of a presheaf
from $\cat{PSM}$ to $\cat{PSM}_{\leq n}$ or to $\cat{PSM}_{n}$;
$\cat{PSM}_{n}$s are nice because of Lemma \ref{lem:PSMn}, while
$\cat{PSM}_{\leq n}$s are nice because they form an increasing chain.

\begin{defn}
Let $u_{n}:\cat{PSM}_{\leq n}\map\cat{PSM}$ be the inclusion. The
\emph{$n$-th approximation} is the functor $\mathrm{app}_{n}=u_{n*}\circ u_{n}^{*}:\widehat{\cat{PSM}}\map\widehat{\cat{PSM}}$.
The \emph{approximating morphism} is the morphism $\mathrm{id}_{\widehat{\cat{PSM}}}\map\mathrm{app}_{n}$
coming from the fact that $u_{n*}$ and $u_{n}^{*}$ form an adjoint
pair of functors. A presheaf $F\in\widehat{\cat{PSM}}$ is \emph{of
order} $n$ if $F\map\mathrm{app}_{n}(F)$ is an isomorphism (it is
then of order $m$ for every $m\geq n$), or equivalently, if it is
of the form $u_{n*}(G)$ for some $G\in\widehat{\cat{PSM}_{\leq n}}$.
\end{defn}
The category of presheaves of order $n$ is thus equivalent to $\widehat{\cat{PSM}_{\leq n}}$
(or $\widehat{\cat{PSM}_{n}}$); in particular, for $n=0$ it is equivalent
to $\widehat{\cat{SM}}$. The chain $\cat{PSM}_{\leq0}\subset\cat{PSM}_{\leq1}\subset\cat{PSM}_{\leq2}\subset\cdots\subset\cat{PSM}$
gives us a chain\[
\mathrm{app}_{0}\leftarrow\mathrm{app}_{1}\leftarrow\mathrm{app}_{2}\leftarrow\cdots\leftarrow\mathrm{app}_{n}\leftarrow\cdots\]
that together with the approximating morphisms $\mathrm{id}_{\widehat{\cat{PSM}}}\map\mathrm{app}_{k}$
form a commutative diagram.

\section{Representability of presheaves}

To get differential worms on an $F$ (see Definition \ref{def:forms-on-functor})
we need $F(\rodd{n},\cdot)$ to be representable, so let us  study
this case.

\begin{defn}
The objects of \emph{the category} $\cat{SM}_{[n]}$ are supermanifolds
with right action of the super-semigroup $\HOM(\rodd{n},\rodd{n})$,
and morphisms are equivariant maps.
\end{defn}
By Lemma \ref{lem:PSMn} this category is a full subcategory of $\widehat{\cat{PSM}_{n}}$,
and thus of $\widehat{\cat{PSM}}$. The objects of $\cat{SM}_{[1]}$
are also called differential non-negatively graded supermanifolds.

The following definition says that $F$ admits level-$n$ worms and
that it is completely determined by these worms.

\begin{defn}
A presheaf $F$ on $\cat{PSM}$ is $n$-\emph{representable} if it
is of order $n$ and if $F(\rodd{n},\cdot)$ is representable.
\end{defn}
$0$-representable presheaves are those represented by objects of
$\cat{PSM}$ of the form $(X,\{\mbox{point}\})$ (where $X$ is the
supermanifold representing $F(\{\mbox{point}\},\cdot)$), i.e.~their
category is equivalent to $\cat{SM}$. Similarly, we have the following
obvious claim for arbitrary $n$:

\begin{prop}
The category of $n$-representable presheaves on $\cat{PSM}$ is equivalent
to $\cat{SM}_{[n]}$, where $Z\in\cat{SM}_{[n]}$ corresponding to
an $n$-representable $F$ is the supermanifold representing $F(\rodd{n},\cdot)$,
and $F$ corresponding to $Z$ is given by \[
F(X,B)=\textrm{the set of }\HOM(\rodd{n},\rodd{n})\textrm{-equivariant maps }\HOM(\rodd{n},X)×B\map Z.\]
In particular the object $Z\in\cat{SM}_{[n]}$ corresponding to $(X,B)$
is $\HOM(\rodd{n},X)×B$. 
\end{prop}
It follows from definition that if $F\in\widehat{\cat{PSM}}$ is of
order $n$ then it is of order $m$ for any $m\geq n$. We shall make
a conjecture that same is valid for $n$-representability, i.e.~that
the categories of $n$-representable functors form an increasing chain.
If the conjecture is valid, by the previous theorem we have a chain
of fully faithful functors\[
\cat{SM}=\cat{SM}_{[0]}\map\cat{SM}_{[1]}\map\cat{SM}_{[2]}\map\cdots\]

\begin{conjecture*}
Let $F$ be an $n$-representable presheaf on $\cat{PSM}$. Then 

1. $F$ is $m$-representable for every $m\geq n$ 

2. $\mathrm{app}_{m}(F)$ is $m$-representable for every $m$.
\end{conjecture*}
\begin{rem*}
1.~follows from 2\@.~since $F\cong\mathrm{app}_{m}(F)$ for $m\geq n$.
We thus have to prove that if $Z$ is an object of $\cat{SM}_{[n]}$
(in our case $Z$ is the supermanifold representing $F(\rodd{n},\cdot)$)
then the following presheaf $H\in\widehat{\cat{SM}}$ is representable:\[
H(X)=\textrm{the set of }\HOM(\rodd{n},\rodd{n})\textrm{-equivariant maps }\HOM(\rodd{n},\rodd{m})×X\map Z.\]
This is easy to see if $m\leq n$, but we were not able to prove (or
disprove) it when $m>n$. (The conjecture is, of course, true for
all the examples we shall consider.)
\end{rem*}
We understood $\widehat{\cat{PSM}}$ and all its subcategories such
as $\widehat{\cat{PSM}_{n}}$ and $\cat{SM}_{[n]}$ as categories
of generalized supermanifolds. Here is the condition for a generalized
\emph{manifold:}

\begin{defn}
\label{def:even-functor}A presheaf $F\in\widehat{\cat{PSM}}$ is
\emph{even} if it is stable under the parity involution, i.e.~if
for every $(X,B)\in\cat{PSM}$, $F$ sends the parity involution $(X,B)\map(X,B)$
(which is the parity involution applied to both $X$ and $B$) to
the identity on $F(X,B)$. 
\end{defn}
If $F$ is represented by some $(Y,C)\in\cat{PSM}$, this condition
means that both $Y$ and $C$ are manifolds. The condition translates
to $\cat{SM}_{[n]}$ in the following way: the parity involution of
$\rodd{n}$ (an element of $\Hom(\rodd{n},\rodd{n})$) should act
on a $Z\in\cat{SM}_{[n]}$ by the parity involution of $Z$. In particular,
for $n=1$ it means that if $f$ is a function on $Z$ of degree $d$,
its parity is the parity of $d$.

\section{Examples of differential forms/worms on presheaves and of approximations\label{sec:Examples}}

The presheaves $F\in\widehat{\cat{PSM}}$ in this section will be
of the form $F(X,B)=$ smooth families of something on $X$, parametrized
by $B$. It's an experimental fact (but not a theorem) that such presheaves
admit worms of any level (i.e.~$F(\rodd{n},\cdot)$ is representable
for every $n$, or in other words, $\mathrm{app}_{n}(F)$ is $n$-representable
for every $n$).

\begin{example*}
As the most trivial example, let us take the presheaf represented
by an object $(X,B)$. This presheaf is $n$-representable for any
$n\geq1$, the corresponding object in $\cat{SM}_{[n]}$ is $\HOM(\rodd{n},X)×B$,
hence the level-$n$ worms on the presheaf are level-$n$ worms on
$X$ smoothly parametrized by $B$. The presheaf is $0$-representable
only if $B$ is a point (which is, however, the most interesting case).
\end{example*}
{}

\begin{example*}
Let $G$ be a Lie group (or supergroup). Let $F(X,B)=\Hom(X×B,G)/\Hom(B,G)$
(where both $\Hom(X×B,G)$ and $\Hom(B,G)$ inherit the group structure
from $G$). Then $F(\rodd{1},\cdot)$ is represented by $\Pi\mathfrak{g}=\HOM(\rodd{1},G)/G$.
Differential forms on $F$, i.e.~the differential graded algebra
of functions on $\Pi\mathfrak{g}$, is the Chevalley-Eilenberg complex
of $\mathfrak{g}$. We have $\mathrm{app}_{1}(F)(X,B)=$ smooth families
of flat $\mathfrak{g}$-connections on $X$ parametrized by $B$;
$\mathrm{app}_{1}(F)$ is thus the sheafification of $F$. All higher
approximations of $F$ are isomorphic to $\mathrm{app}_{1}(F)$.
\end{example*}
{}

\begin{example*}
As a generalization of the previous example, let $G$ be a Lie groupoid
(or supergroupoid). From $X,B\in\cat{SM}$ we form a groupoid $G_{X,B}$
with objects $X×B$ and morphisms $X×X×B$ (i.e.~the pair groupoid
of $X$ times the trivial groupoid of $B$). We define \[
F(X,B)=\{\textrm{Lie groupoid morphisms }G_{X,B}\map G\}\]
(in other words, families of Lie groupoid morphisms from the pair
groupoid of $X$ smoothly parametrized by $B$). Then $F(\rodd{1},\cdot)$
is represented by $\Pi A$, where $A$ is the Lie algebroid of $G$.
The algebra of differential forms on $F$ is thus $\Gamma(\bigwedge A)$,
the Chevalley-Eilenberg complex of $A$. The first approximation of
$F$ is given by \[
\mathrm{app}_{1}(F)(X,B)=\{\textrm{Lie algebroid morphisms }TX×B\map G\},\]
i.e.~it is the sheafification of $F$. Higher approximations of $F$
are again isomorphic to $\mathrm{app}_{1}(F)$.
\end{example*}
{}

\begin{example*}
This is an extremely simple example, but it might be enlightening.
Let $F$ be given by \[
F(X,B)=\{\textrm{families of closed }k\textrm{-forms on }X\textrm{ smootly parametrized by }B\}.\]
Then $F(\rodd{1},\cdot)$ is represented by the (1-dimensional) vector
superspace of closed $k$-forms on $\rodd{1}$, i.e.~by $\mathbb{R}[k]$.
This $F$ is $1$-representable.
\end{example*}
{}

\begin{example*}
This is an example of a presheaf on which a group acts. Let $G$ be
a Lie group, $\mathfrak{g}$ its Lie algebra, and let $F$ be given
by\[
F(X,B)=\{\textrm{families of }\mathfrak{g}\textrm{-connections on }X\textrm{ smoothly parametrized by }B\}.\]
 On this $F$ acts the group (represented by) $G$, namely on $F(X,B)$
acts the group $\Hom(X\times B,G)$ by gauge transformations. The
presheaf $F$ is $1$-representable, the algebra of differential forms
on $F$ is the Weil algebra $W(\mathfrak{g})$. Moreover, the action
of $\HOM(\rodd{1},G)$ on differential forms gives the standard $G$-differential
structure on $W(\mathfrak{g})$.

If now $M$ is a manifold (or supermanifold) with $G$-action, the
algebra of differential forms on $F\times M$ is $W(\mathfrak{g})\otimes\Omega(M)$,
and the algebra of $\HOM(\rodd{1},G)$-invariant differential forms
is the basic subcomplex (whose cohomology is the equivariant cohomology
of $M$ if $G$ is compact).
\end{example*}
{}

\begin{example*}
Let us fix a supermanifold $Y$ and let $F$ be given by $F(X,B)=\Hom(X×X×B,Y)$
(or in other words, families of smooth maps $X\times X\map Y$ smoothly
parametrized by $B$). Then $F(\rodd{n},\cdot)$ is represented by
$\HOM(\rodd{n}\times\rodd{n},Y)$, and \[
\mathrm{app}_{n}(F)(X,B)=\{\textrm{smooth families of sections of }j^{k}(X,Y)\textrm{ parametrized by }B\},\]
where $j^{k}(X,Y)\map X$ is the bundle of $k$-jets of maps $X\map Y$.
If for example $Y=\mathbb{R}$, the differential forms on $F$ form
the free graded commutative algebra generated by $x,\xi,\tau,t$ with
$\deg x=0$, $\deg\xi=\deg\tau=1$, $\deg t=2$, with differential
given by $dx=\xi$, $d\tau=t$.
\end{example*}
{}

\begin{example*}
Given $m\in\mathbb{N}$ we shall consider the presheaf given by

\[
F(X,B)=\{\textrm{smooth families of sections of }(T^{*})^{\otimes m}X\textrm{, parametrized by }B\},\]
and similar presheaves given by natural subbundles of $(T^{*})^{\otimes m}X$
(it is convenient to decompose $(T^{*})^{\otimes m}X$ with respect
to the action of the symmetric group $S_{m}$).

Let $\lambda$ be a Young diagram with $m$ squares, e.g. 
$$\includegraphics{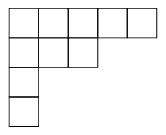}$$
for $m=10$, and let $W_{\left(\lambda\right)}$ be the corresponding
irreducible representation of $S_{m}$. If $V$ is a vector space,
let \[
V_{\lambda}=\Hom_{S_{m}}(W_{\left(\lambda\right)},V^{\otimes m}),\]
where $S_{m}$ acts on $V^{\otimes m}$ by permutations of factors
(recall that $V_{\lambda}$ is an irreducible representation of $GL(V)$).
This way we get a functor from the category of vector spaces to itself,
given by $V\mapsto V_{\lambda}$, and thus a presheaf \[
F_{\lambda}(X,B)=\{\textrm{smooth families of sections of }T_{\lambda}^{*}X\textrm{, parametrized by }B\}.\]

For any $\lambda$ and any $n$ the functor $F_{\lambda}(\rodd{n},\cdot)$
is represented by the vector superspace $\Gamma(T_{\lambda}^{*}\rodd{n})$.
Let $c$ be the number of columns of $\lambda$. The presheaf $F_{\lambda}$
is $n$-representable for $n>c$, while $\mathrm{app}_{n}(F_{\lambda})(X,B)=\{0\}$
for $n<c$. For $n=c$ we have \[
\mathrm{app}_{c}(F_{\lambda})(X,B)=\{\textrm{smooth families of sections of }\widetilde{T_{\lambda}^{*}}X\textrm{, parametrized by }B\},\]
where $\widetilde{T_{\lambda}^{*}}X$ is a natural bundle over $X$
containing $T_{\lambda}^{*}X$.

If $c=1$ (the case of differential forms) then $\widetilde{T_{\lambda}^{*}}X=T_{\lambda}^{*}X$;
the same is true when $\lambda$ has only one row (the case of symmetric
tensors). When $c=2$ we have the following result: $\widetilde{T_{\lambda}^{*}}X$
has a natural increasing filtration\[
T_{\lambda}^{*}X=A_{1}\subset A_{2}\subset\cdots\subset\widetilde{T_{\lambda}^{*}}X,\]
such that each $A_{i+1}/A_{i}$ is isomorphic to some $T_{\mu}^{*}X$.
The rule for obtaining these $\mu$'s from $\lambda$ should be clear
from this picture:$$
\includegraphics{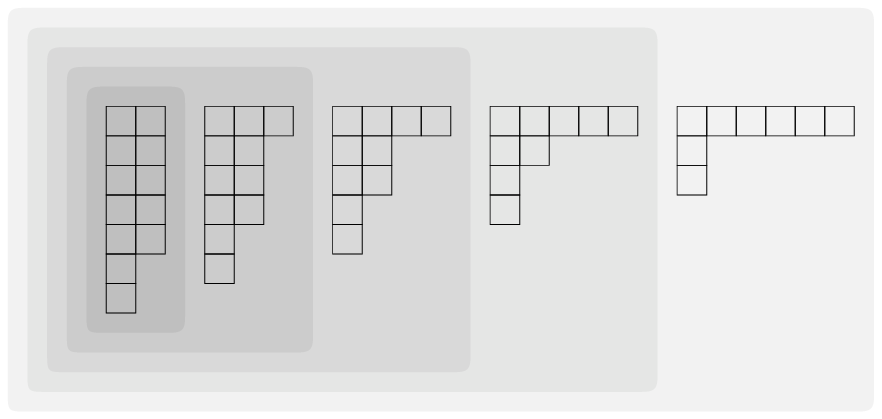}
$$

In words, we keep removing squares from the two column and adding
them to the first row, till the second column contains only one square.

This result enables us to decompose $\Omega_{[2]}(M)$ to indecomposable
representation of the super-semigroup $\HOM(\rodd{2},\rodd{2})$.
The irreducible (left) representations of $\HOM(\rodd{2},\rodd{2})$
are duals of $\Gamma(T_{\lambda}^{*}\rodd{2})$ for 2-column $\lambda$'s
(these are called generic irreducibles), and duals of the  spaces
of closed $k$-forms on $\rodd{2}$, $k\in\mathbb{N}$ (non-generic
irreducibles). The representation theory of $\HOM(\rodd{2},\rodd{2})$
is quite simple: generic irreducibles never appear in the composition
series of reducible indecomposable representations, and the only reducible
indecomposable cyclic representations are duals to the spaces of differential
$k$-forms on $\rodd{2}$, $k\in\mathbb{N}$. It gives the following
decomposition of polynomial functions on $\HOM(\rodd{2},\rodd{2})$:\[
\Omega_{[2]}(\rodd{2})\cong\left(\bigoplus_{\lambda}\Gamma(T_{\lambda}^{*}\rodd{2})^{*}\otimes\Gamma(T_{\lambda}^{*}\rodd{2})\right)\oplus\frac{\bigoplus_{k}\Omega^{k}(\rodd{2})^{*}\otimes\Omega^{k}(\rodd{2})}{d(\bigoplus_{k}\Omega^{k}(\rodd{2})^{*}\otimes\Omega^{k-1}(\rodd{2}))}\]
(where the sum is over all 2-column $\lambda$'s) and thus for any
supermanifold $X$\[
\Omega_{[2]}(X)\cong\left(\bigoplus_{\lambda}\Gamma(T_{\lambda}^{*}\rodd{2})^{*}\otimes\Gamma(\widetilde{T_{\lambda}^{*}}X)\right)\oplus\frac{\bigoplus_{k}\Omega^{k}(\rodd{2})^{*}\otimes\Omega^{k}(X)}{d(\bigoplus_{k}\Omega^{k}(\rodd{2})^{*}\otimes\Omega^{k-1}(X))}.\]

Unfortunately, the representation theory of $\HOM(\rodd{n},\rodd{n})$
is wild for $n\geq3$ \cite{sh}, and we do not know how to decompose
the space of polynomial functions on this semigroup. We also do not
know the composition series of $\widetilde{T_{\lambda}^{*}}$ for
general $\lambda$ with more than 2 columns. Just as an example, here
is the composition series for a rather simple $\lambda$ with 3 columns:
$$
\includegraphics{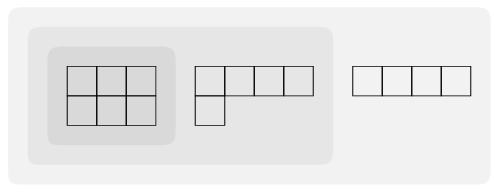}
$$
\end{example*}

\section{The case of sheaves}

Most of the examples of presheaves we considered were actually sheaves,
and all their approximations were sheaves. Let us be more specific
about what we mean. Let $(X,B)$ and $(X_{i},B_{i})$ (where $i$
runs through some index set $I$) be objects of $\mathsf{PSM}$. Then
$(X_{i},B_{i})$ is a \emph{covering} of $(X,B)$ if each $X_{i}$
is an open subset of $X$, each $B_{i}$ is an open subset of $B$
and $\bigcup(X_{i}\times B_{i})=X\times B$. A presheaf on $\mathsf{PSM}$
is a \emph{sheaf} if for any $(X,B)$ and any covering of $(X,B)$
it satisfies the usual conditions for sheaves.

Sheaves on $\mathsf{PSM}$ can be extended to a larger category where
they look more natural:

\begin{defn}
The category $\cat{FSM}$ is the category of foliated supermanifolds,
i.e.~its objects are foliated supermanifolds and morphisms are smooth
maps that map leaves into leaves.
\end{defn}
Any object $(X,B)\in\cat{PSM}$ gives a foliated supermanifold $X\times B$,
with leaves given by projection $X\times B\map B$. This way $\cat{PSM}$
becomes a full subcategory of $\cat{FSM}$. Any sheaf on $\cat{PSM}$
is naturally extended to a presheaf on $\cat{FSM}$ that is actually
a sheaf w.r.t.~the \'etale topology.

As an important example, all $n$-representable presheaves (for any
$n$) are sheaves. All approximations of sheaves are sheaves as well.

\section{Open ends}

This concluding section contains some suggestions for further development.

\subsection{Stacks, their differential forms/worms and approximations}

Very informally speaking, the basic idea of this paper is to approximate
{}``things'' defined on all (super)manifolds by evaluating them
on the simplest supermanifolds of the form $\rodd{k}$; we suppose
these {}``things'' to be contravariant, thus getting differential
forms/worms into the play. There are of course many different formalizations
of this vague idea.

As an example, here we shall briefly present a slight generalization
of what we did so far, basically by passing from contravariant functors
to the category of sets to contravariant functors to the category
of categories. More precisely, we shall deal with stacks (sheaves
with values in categories) over the category $\cat{SM}$ or $\cat{FSM}$
(we always use the étale topology, where by definition a covering
is a surjective locally-diffeomorphic map). We do it just to show
some interesting examples. For definitions see \cite{gir,sga1}, though
they are not necessary to understand our simple examples.

If $G$ is a Lie group, we have the stack of all principal $G$-bundles
(morphisms in this stack are $G$-equivariant maps between the bundles).
More generally, if $G$ is a Lie supergroup with a right action of
$\HOM(\rodd{k},\rodd{k})$ (i.e.~$G$ is a group object in the category
$\cat{SM}_{[k]}$) then for any (super)manifold $X$ we can consider
principal $G$-bundles \[
\begin{array}{c}
P\\
\downarrow\\
\HOM(\rodd{k},X)\end{array}\]
 with a compatible action of $\HOM(\rodd{k},\rodd{k})$ on $P$ (compatibility
means that the maps $G\times P\map P$ and $P\map\HOM(\rodd{k},X)$
are $\HOM(\rodd{k},\rodd{k})$-equivariant). We again get a stack
over $\cat{SM}$ (morphisms are now maps that are both $G$ and $\HOM(\rodd{k},\rodd{k})$
equivariant). We can extend this stack to $\cat{FSM}$ by using the
same formulas, when $\rodd{k}$ is understood as a foliated supermanifold
with just one leaf. 

Instead of a Lie group $G$ we can also take any Lie category $C$
and consider principal $C$-bundles. If $C\in\cat{SM}_{[n]}$ (the
right $\HOM(\rodd{k},\rodd{k})$-action on $C$ is supposed to be
copatible with the category structure of $C$), we can again consider
principal $C$-bundles over $\HOM(\rodd{k},X)$'s with a compatible
right $\HOM(\rodd{k},\rodd{k})$-action, and get a stack over $\cat{FSM}$
this way. 

Yet more generally, the action of $\HOM(\rodd{k},\rodd{k})$ on $C$
may be weak (i.e.~up to coherent natural transformations of $C$);
precisely and invariantly it means that we have a fibred Lie category
$\hat{C}\map\HOM(\rodd{k},\rodd{k})$ (where the semigroup $\HOM(\rodd{k},\rodd{k})$
is understood as a category with only one object) whose fibre is $C$.
We can then consider principal $C$-bundles $P\map\HOM(\rodd{k},X)$
with a compatible action of $\hat{C}$ on $P$. Stacks of this form
will be called \emph{weakly $k$-representable}. 

Contrary to $k$-representable (pre)sheaves, the fibred category $\hat{C}$
is not uniquely determined by the corresponding stack on $\cat{FSM}$.
However, given a morphism between two stacks represented by $\hat{C_{1}}$
and $\hat{C_{2}}$ respectively, we have a Lie category fibred over
$(1\map2)\times\HOM(\rodd{k},\rodd{k})$ (where $(1\map2)$ is the
category with objects 0 and 1, with a unique morphism $1\map2$) that
coincides with $\hat{C_{1}}$ and $\hat{C_{2}}$ over 1 and 2 respectively.
If the morphism is an equivalence, we get a Lie category fibred over
$(1\leftrightarrow2)\times\HOM(\rodd{k},\rodd{k})$; the category
$\hat{C}$ is defined by the corresponding stack up to such equivalences.

We can define approximations of stacks on $\cat{FSM}$ just as we
did for presheaves: given a stack $\cat{S}\map\cat{FSM}$ we can restrict
it to $\cat{PSM}_{k}$; the restriction has a right adjoint (see \cite{gir}).
The composition of the restriction with the induction will be called
again $k$-th approximation. Weakly $k$-representable stacks are
of course equivalent to their $k$-th (and higher) approximations.

Let us finally present some examples. They are all weakly 1-representable.

\begin{example*}
\textbf{({}``Categorified de Rham complex'')} Let $\mathbb{R}[k]=\mathcal{Z}^{k}(\rodd{1})\subset\Omega^{k}(\rodd{1})$
denote the 1-dimensional group of closed $k$-forms on $\rodd{1}$.
If $X$ is any supermanifold, it is easy to see that equivariant principal
$\mathbb{R}[k]$-bundles $P\map\HOM(\rodd{1},X)$ are classified by
$k+1$-th de Rham cohomology group of $X$, while automorphisms of
$P$ are canonically isomorphic to closed $k$-forms on $X$ (see
\cite{s1}). These bundles are nothing but $\mathcal{Z}^{k}(X)$-torsors.
(If $X\in\cat{FSM}$ we need to take leafwise forms and their cohomology.)

On the other hand, equivariant $\Omega^{k}(\rodd{1})$-principal bundles
$P\map\HOM(\rodd{1},X)$ are the same as $\Omega^{k}(X)$-torsors,
i.e.~affine bundles over $X$ whose associated vector bundle is $\bigwedge^{k}T^{*}X$.
They are all isomorphic.

The exact sequence\[
0\map\mathbb{R}[k]\map\Omega^{k}(\rodd{1})\map\mathbb{R}[k+1]\map0\]
says that a reduction of a $\Omega^{k}(\rodd{1})$-bundle $P\map\Hom(\rodd{1},X)$
(i.e.~of a $\Omega^{k}(X)$-torsor) to a $\mathbb{R}[k]$-bundle
(to a $\mathcal{Z}^{k}(X)$-torsor) is the same as a trivialization
of the $\mathbb{R}[k+1]$-bundle $X/\mathbb{R}[k]$.
\end{example*}
{}

\begin{example*}
\textbf{(Pontryagin class)} Let $G$ be a Lie group and $\left\langle ,\right\rangle $
be an invariant symmetric bilinear form on the Lie algebra $\mathfrak{g}$.
There is a central extension $\widetilde{T[1]}\mathfrak{g}$ of the
cone $T[1]\mathfrak{g}=\mathfrak{g}\oplus\mathfrak{g}[1]$ of $\mathfrak{g}$
by $\mathbb{R}[2]$, given by \[
[u,v]=\left\langle u,v\right\rangle \in\mathbb{R}[2],\quad u,v\in\mathfrak{g}[1].\]
This extention of differential graded Lie algebras can then be integrated
to a central extension $\widetilde{T[1]}G$ of $T[1]G$ by $\mathbb{R}[2]$.

Since on both $\widetilde{T[1]}\mathfrak{g}$ and $\widetilde{T[1]}G$
we have a right action of $\HOM(\rodd{1},\rodd{1})$, they give us
1-representable sheaves on $\cat{FSM}$: for any foliated supermanifold
$X$, $\widetilde{T[1]}\mathfrak{g}$ gives a central extension of
the Lie algebra of all smooth maps $X\map\mathfrak{g}$ by closed
leafwise 2-forms on $X$, while $\widetilde{T[1]}G$ gives a central
extension of the group of all smooth maps $X\map G$, again by closed
leafwise 2-forms on $X$.

The stack weakly 1-represented by $\widetilde{T[1]}G$ is the stack
of torsors of the above-mentioned sheaf of group. If $P\map X$ is
a principal $G$-bundle, it can be lifted to a torsor of the central
extension iff its $\left\langle ,\right\rangle $-Pontryagin class
vanishes, i.e.~if the leafwise 4-form $\left\langle \Omega\wedge\Omega\right\rangle $
(where $\Omega$ is the curvature of a leafwise connection on $P$)
is exact.
\end{example*}
{}

\begin{example*}
\textbf{(Weil and Cartan models of equivariant cohomology)} Let $G$
be a Lie group. Let us consider the stack whose objects are principal
$G$-bundles with connection (if we mean a stack over $\cat{FSM}$,
the connection should be just along the leaves of the foliation) and
morphisms are $G$-equivariant maps preserving the connection. This
stack is weakly 1-representable. It is easy to find a Lie supergoupoid
$\Gamma$ with a right action of $\HOM(\rodd{1},\rodd{1})$ representing
this stack: over $\rodd{1}$ any $G$-bundle is trivial, so we can
take the affine space $\Omega^{1}(\rodd{1})\otimes\mathfrak{g}$ of
connections on this bundle, the gauge group $\HOM(\rodd{1},G)$ acting
on this space, and finally put $\Gamma$ to be the corresponding action
groupoid.

We can find a smaller groupoid $\Gamma'$ equivalent to $\Gamma$:
any $\mathfrak{g}$-connection on $\rodd{1}$ can be made to vanish
at the origin by a gauge transformation, i.e.~is gauge-equivalent
to a connection of the form $a\,\theta d\theta$, $a\in\mathfrak{g}$.
The groupoid $\Gamma'$ will be the full subgroupoid of $\Gamma$
whose objects are connections of this form. $\Gamma'$ is nothing
but the action groupoid of $G$ on $\mathfrak{g}$ (with the adjoint
action). The action of $\HOM(\rodd{1},\rodd{1})$ doesn't descend
from $\Gamma$ to $\Gamma'$ however; we only have a fibred Lie category
$\hat{\Gamma'}\map\HOM(\rodd{1},\rodd{1})$ with fibre $\Gamma$.

Let now $M$ be a manifold with a $G$-action. Let us consider the
stack of principal $G$ bundles with a connection (as above) and with
a $G$-equivariant map to $M$. This stack is weakly 1-represented
by the action groupoid of $\HOM(\rodd{1},G)$ acting on $\Omega^{1}(\rodd{1})\otimes\mathfrak{g}\times\HOM(\rodd{1},M)$,
or by the equivalent action groupoid of $G$ on $\mathfrak{g}\times\HOM(\rodd{1},M)$.
When we take the invariant functions on the bases of these groupoids,
we get Weil and Cartan model of $G$-equivariant cohomology of $M$
respectively. The fact that in Cartan model $d^{2}=0$ only on equivariant
forms (i.e.~on invariant functions on $\mathfrak{g}\times\HOM(\rodd{1},M)$)
comes from the fact that on $\Gamma'$ we do not have a $\HOM(\rodd{1},\rodd{1})$-action,
just the fibred category $\hat{\Gamma'}\map\HOM(\rodd{1},\rodd{1})$.
\end{example*}
{}

\begin{example*}
\textbf{(Quasi-Poisson groupoids)} An interesting example of Lie categories
fibred over $\HOM(\rodd{1},\rodd{1})$, and thus of weakly 1-representable
stacks, comes from Lie quasi-bialgebras and more generally from Lie
quasi-bialgebroids. As at the end of the previous example, we will
\emph{not} have an action of $\HOM(\rodd{1},\rodd{1})$ on a Lie category,
just a Lie category fibred over $\HOM(\rodd{1},\rodd{1})$. More precisely,
we shall use the isomorphism $\HOM(\rodd{1},\rodd{1})=(\mathbb{R},\times)\ltimes\mathbb{R}[-1]$
(a general element of $\HOM(\rodd{1},\rodd{1})$ is of the form $\theta\mapsto a\theta+\beta$;
$(\mathbb{R},\times)$ corresponds to $\theta\mapsto a\theta$ and
$\mathbb{R}[-1]$ to $\theta\mapsto\theta+\beta$); we shall describe
Lie groupoids fibred over $\mathbb{R}[-1]$, with a right action of
$(\mathbb{R},\times)$ (in other words, we shall have graded Lie groupoids,
with differentials {}``up to homotopy'').

As noticed in \cite{s2}, a Lie quasi-bialgebra structure on a vector
space $\mathfrak{g}$ can be described as a graded principal $\mathbb{R}[2]$-bundle
$X\map\mathfrak{g}^{*}[1]$, with an $\mathbb{R}[2]$-invariant odd
Poisson structure $\pi$ of degree $-1$ on $X$. Similarly, if $A\map M$
is a vector bundle, a Lie quasi-bialgebroid structure on $A$ is a
principal $\mathbb{R}[2]$-bundle $X\map A^{*}[1]$, again with an
odd Poisson structure as above. (In perhaps more familiar algebraic
terms, we add a variable $t$ of degree 2 to the graded-commutative
algebra $\Gamma(\bigwedge A^{*})$ and want a Gerstenhaber bracket
on the result, such that all $[\alpha,\beta]$, $[\alpha,t]$ and
$[t,t]$ are in $\Gamma(\bigwedge A^{*})$ for any $\alpha,\beta\in\Gamma(\bigwedge A^{*})$.)
If the odd Poisson structure on $X$ is integrable to an odd symplectic
groupoid $\hat{\Gamma}$, the $\mathbb{R}[2]$-action on $X$ lifts
to a Hamiltonian action on $\hat{\Gamma}$, with a Hamiltonian $\hat{\Gamma}\map\mathbb{R}[-1]$
that is a morphism of groupoids and  makes $\hat{\Gamma}$ to a groupoid
fibred over $\mathbb{R}[-1]$, as we wanted.
\end{example*}
{}

\subsection{Differences vs.~differentials}

In this paper we exploited the fact that $\HOM(\rodd{k},X$) is representable.
Besides $\rodd{k}$'s there are other (super)manifolds sharing this
property, namely finite sets (of course, for them it is completely
trivial). We could thus rewrite the paper: in place of $\rodd{k}$
we would use a $m$-element set, in place of $\HOM(\rodd{k},\rodd{k})$
we would use the semigroup $\Sigma_{m}$ of all maps of the $m$-element
set to itself etc.

In the Examples Section \ref{sec:Examples} some of the examples are
induced from finite sets. Let us consider the second one given by
a Lie group $G$ (where we can take $m=3$). When we restrict that
functor to the category $\Sigma$ of finite sets we get a functor
$\Sigma^{o}\map\cat{SM}$ (for any finite set $S$ we take its pair
groupoid $S\times S\rightrightarrows S$ and then the manifold of
all the functors from this groupoid to $G$). The functor $\Sigma^{o}\map\cat{SM}$
is roughly speaking the {}``enhanced'' nerve of $G$ (the ordinary
nerve is a contravariant functor from the category of ordered finite
sets $\Delta$; in $\Sigma$ we have more morphisms than in $\Delta$,
i.e.~a functor from $\Sigma$ is a stronger structure than a functor
from $\Delta$; we get {}``simplicial objects with inverses'').
The general procedure of computing {}``the Lie algebra of a simplicial
manifold with inverses'' is to start with a functor $\Sigma^{o}\map\cat{SM}$,
induce it to a presheaf on $\cat{PSM}$ and then compute its differential
forms and possibly also its higher level worms. An interesting question
is when these algebras exist (i.e.~when the corresponding presheaves
on $\cat{SM}$ are representable) and whether/when the $k$-th approximations
of this presheaf stabilize.

\appendix

\section{Cohomology of worms}

The cohomology of $\Omega_{[k]}(X)$ with respect to $d_{1}$ (or
with respect to any non-zero linear combination of $d_{a}$'s) is
canonically isomorphic to de Rham cohomology of $X$. Indeed, $d_{1}$
comes from the vector field $\partial/\partial\theta_{1}$ on $\rodd{1}$.
On the other hand, \[
[\partial/\partial\theta_{1},\theta_{1}\theta_{2}\,\partial/\partial\theta_{2}]=\theta_{2}\ \partial/\partial\theta_{2},\]
which means that the $d_{1}$-cohomology of the subcomplex of $\Omega_{[k]}(X)$
with non-zero 2nd degree vanishes (as $\theta_{2}\ \partial/\partial\theta_{2}$
generates the 2nd degree, and the above equation says that it is homotopic
to 0). The same is of course true for any degree except for the first
one. As a result, the embedding $\Omega(X)=\Omega_{[1]}(X)\map\Omega_{[k]}(X)$
is a quasiisomorphism.

The action of $\HOM(\rodd{k},\rodd{k})$ gives rise also to some other
differentials on $\Omega_{[k]}(X)$. For example, the cohomology of
the differential given by $\theta_{1}\theta_{2}\ \partial/\partial\theta_{2}$
is isomorphic to $\Omega_{[k-1]}(X)$: the above equation again shows
that $\Omega_{[k]}(X)$ is quasiisomorphic to the subcomplex with
vanishing 2nd degree (which can be identified with $\Omega_{[k-1]}(X)$),
but on the subcomplex the differential vanishes.

\section{Cartan wormulas}

Here we shall see how the Cartan formula\[
[d,i_{v}]=\mathcal{L}_{v}\]
can be derived from the fact that differential forms are functions
on a map space, and how the formula generalizes to $\Omega_{[k]}$.
The formula really comes from the action of the group $\mathit{Diff}(Y)\ltimes(\mathit{Diff}(X))^{Y}$
on the space of maps $Y\map X$ (where $Y=\rodd{k}$ for our purposes),
or yet better, from the action of the corresponding Lie algebra $\mathcal{X}(Y)\ltimes(C^{\infty}(Y)\otimes\mathcal{X}(X))$
(here $\mathcal{X}(X)$ denotes the Lie algebra of vector fields on
$X$).

If $u\in\mathcal{X}(Y)$ and $f\otimes v\in C^{\infty}(Y)\otimes\mathcal{X}(X)$,
we denote the corresponding vector fields on $\HOM(Y,X)$ by $u^{\flat}$
and $f\cdot v^{\#}$ ($1\cdot v^{\#}$ will be denoted simply $v^{\#}$).%
\footnote{one can easily see that the operation $f\cdot$ can be applied to
any vector field on $\HOM(Y,X)$, not just to those of the form $v^{\#}$;
we shall not need this fact here%
} The following formulas express the fact that on $\HOM(Y,X)$ we have
an action of $\mathcal{X}(Y)\ltimes(C^{\infty}(Y)\otimes\mathcal{X}(X))$:\begin{eqnarray*}
[u_{1}^{\flat},u_{2}^{\flat}] & = & [u_{1},u_{2}]^{\flat}\\
{}[u^{\flat},f\cdot v^{\#}] & = & (uf)\cdot v^{\#}\\
{}[f_{1}\cdot v_{1}^{\#},f_{2}\cdot v_{2}^{\#}] & = & (-1)^{|v_{1}||f_{2}|}(f_{1}f_{2})\cdot[v_{1},v_{2}]^{\#}.\end{eqnarray*}

In the case of $Y=\rodd{1}$ the equation \[
[(\partial/\partial\theta)^{\flat},\theta\cdot v^{\#}]=v^{\#}\]
 is Cartan's $[d,i_{v}]=\mathcal{L}_{v}$ (since $(\partial/\partial\theta)^{\flat}=-d$,
$\theta\cdot v^{\#}=-i_{v}$ and $v^{\#}=\mathcal{L}_{v}$); for $Y=\rodd{k}$,
the above commutation relations are the promissed generalization of
Cartan's formula.

\section{Integration of worms}

Integration of forms can be reformulated in the following well-known
way: For any supermanifold $X$ there is a natural volume form $\mu_{\Pi TX}$
on $\Pi TX$. If $x^{i}$'s are local coordinates on $X$ then $\mu_{\Pi TX}$
is just the coordinate volume form in the coordinates $x^{i}$'s,
$dx^{i}$'s on $\Pi TX$. The volume form is $d$-invariant, i.e.~$\mathcal{L}_{d}\ \mu_{\Pi TX}=0$.
If $\alpha$ is a (pseudo)differential form on $X$, i.e.~a function
on $\Pi TX$, then the integral of $\alpha$ is defined as \[
\int_{\Pi TX}\alpha\ \mu_{\Pi TX}.\]
From the $d$-invariance of $\mu_{\Pi TX}$ we get \[
\int d\alpha=0\]
for any $\alpha$ with compact support. To get Stokes theorem from
here, let $\chi_{\Omega}$ be the characteristic function of a compact
domain $\Omega$; then \[
0=\int d(\chi_{\Omega}\alpha)=\int(d\chi_{\Omega})\alpha+\int\chi_{\Omega}d\alpha=-\int_{\partial\Omega}\alpha+\int_{\Omega}d\alpha.\]

For $k\geq2$ we can write $(\Pi T)^{k}X$ as $\Pi T((\Pi T)^{k-1}X)$
and thus obtain a volume form $\mu_{(\Pi T)^{k}X}$ on $(\Pi T)^{k}X$.
A simple computation shows that $\mu_{(\Pi T)^{k}X}$ is $\mathcal{X}(\rodd{k})$
invariant (only if $k\geq2$; for $k=1$ we had just $d$-invariance,
but not $\theta\,\partial/\partial\theta$-invariance). We use this
volume form to define integrals of worms. For example, if $X=\mathbb{R}$
and $k=2$, we have\[
\int e^{-x^{2}-(d_{1}d_{2}x)^{2}}d_{1}x\, d_{2}x=\pi.\]
The $\mathcal{X}(\rodd{k})$-invariance leads to obvious generalizations
of Stokes theorem. Notice that the elements of $\Omega_{[k]}(M)$
for $k\geq2$ are never integrable functions on $(\Pi T)^{k}M$ (being
polynomial on the fibres of $(\Pi T)^{k}M$). To get interesting examples
we need to use general (non-polynomial) functions on $(\Pi T)^{k}M=\HOM(\rodd{k},M)$;
we shall call them \emph{pseudodifferential worms}.

Here is a simple illustrative result for $k=2$:

\begin{prop}
\label{pro:integ-char}Let $M$ be a connected compact manifold and
$\beta$ a level-2 pseudodifferential worm on $M$. If $d_{1}\beta=d_{2}\beta=0$
and if $\beta$ is integrable then\[
\int\beta=\frac{\beta|_{M}}{2}(-\pi)^{m/2}S_{m}\,\chi_{M},\]
where $m$ is the dimension of $M$, $\beta|_{M}$ is the restriction
of $\beta$ to the zero-section $M\subset(\Pi T)^{2}M$ ($\beta|_{M}$
is a constant because of $d_{1}\beta=d_{2}\beta=0$), $\chi_{M}$
is the Euler characteristics of $M$ and $S_{m}$ is the volume of
the unit $m$-dimensional sphere.
\end{prop}
\begin{proof}
We will first prove a special case. Let $g$ be a Riemann metric on
$M$ and let us define $\gamma=g_{ij}\, d_{1}x^{i}d_{2}x^{j}\in\Omega_{[2]}(M)$;
$\gamma$ is clearly independent of the choice of local coordinates.
We will prove the theorem for $\beta=e^{d_{1}d_{2}\gamma}$. In Riemann
normal coordinates we have at the origin\[
d_{1}d_{2}\gamma=-\delta_{ij}\, d_{1}d_{2}x^{i}d_{1}d_{2}x^{j}-\frac{1}{2}R_{ijkl}\, d_{1}x^{i}\, d_{1}x^{j}\, d_{2}x^{k}\, d_{2}x^{k}\]
where $R_{ijkl}$ are the components of the curvature tensor. To compute
the integral we pass to Riemann normal coordinates at any point of
$M$ and then integrate over $d_{1}x^{i}$'s, $d_{2}x^{i}$'s and
$d_{1}d_{2}x^{i}$'s; we end up with the Pfaffian of the curvature
whose integral over $M$ is well known to be a multiple of $\chi_{M}$.
By comparing with the case when $M$ is the unit $m$-dimensional
sphere we get our result.

To prove the theorem for a general $\beta$ we need the following
result: if $\alpha\in C^{\infty}((\Pi T)^{2}M)$ grows at most polynomially
in $d_{1}d_{2}x^{i}$'s, $d_{1}\alpha=d_{2}\alpha=0$ and $\alpha|_{M}=0$
then $\int e^{d_{1}d_{2}\gamma}\alpha=0$. Indeed, since $\alpha|_{M}=0$,
we can find a $\kappa\in C^{\infty}((\Pi T)^{2}M)$ such that $\alpha=E\kappa$
and $d_{1}\kappa=d_{2}\kappa=0$, where the vector field $E=-(\theta_{1}\partial/\partial\theta_{1}+\theta_{2}\partial/\partial\theta_{2})^{\flat}$
is the total degree, i.e. $E=d_{a}x^{i}\,\partial/\partial(d_{a}x^{i})+2\, d_{1}d_{2}x^{i}\,\partial/\partial(d_{1}d_{2}x^{i})$.
Since \[
\theta_{1}\partial/\partial\theta_{1}+\theta_{2}\partial/\partial\theta_{2}=[\partial/\partial\theta_{1},\theta_{1}\theta_{2}\,\partial/\partial\theta_{1}]+[\partial/\partial\theta_{2},\theta_{2}\theta_{1}\,\partial/\partial\theta_{2}],\]
i.e.\[
E=[d_{1},(\theta_{1}\theta_{2}\,\partial/\partial\theta_{1})^{\flat}]+[d_{2},(\theta_{2}\theta_{2}\,\partial/\partial\theta_{2})^{\flat}],\]
we have \[
\alpha=E\kappa=d_{1}(\theta_{1}\theta_{2}\,\partial/\partial\theta_{1})^{\flat}\kappa+d_{2}(\theta_{2}\theta_{2}\,\partial/\partial\theta_{2})^{\flat}\kappa,\]
and thus\[
\int e^{d_{1}d_{2}\gamma}\alpha=\int d_{1}(e^{d_{1}d_{2}\gamma}(\theta_{1}\theta_{2}\,\partial/\partial\theta_{1})^{\flat}\kappa)+\int d_{2}(e^{d_{1}d_{2}\gamma}(\theta_{2}\theta_{1}\,\partial/\partial\theta_{2})^{\flat}\kappa)=0\]
(the last equality by the $d_{1,2}$-invariance of integration).

Finally to prove our theorem, set $\alpha=\beta-\beta|_{M}$, multiply
$\gamma$ by a constant $s>0$ and take the limit $s\map0_{+}$.
\end{proof}

\section{Level-2 worms in Riemannian geometry}

This appendix is a very elementary example of how worms may appear
in traditional geometry. At the same time it contains an explicit
calculation of an approximation of a presheaf in a simple case, hopefully
better explaining the meaning of approximations than the main abstract
text. We try to be {}``pedagogical'' here.

Given a Riemann metric $g$ on $M$ we can form a level-2 worm \[
\gamma=g_{ij}\, d_{1}x^{i}d_{2}x^{j}\in\Omega_{[2]}(M)\]
 ($\gamma$ is clearly independent of the choice of local coordinates;
we used this worm already in the proof of Proposition \ref{pro:integ-char}).
Out of $\gamma$ we can compute $d_{1}\gamma$, $d_{2}\gamma$ and
$d_{1}d_{2}\gamma$; the last one contains {}``conveniently packed''
Levi-Civita connection and curvature of $g$. Moreover, as we'll see,
the four worms form a basis of an irreducible representation of $\HOM(\rodd{2},\rodd{2})$.

An easy computation gives\begin{equation}
d_{1}d_{2}\gamma=-g_{ij}\, d_{1}d_{2}x^{i}\, d_{1}d_{2}x^{j}+2\,\Gamma_{ijk}\, d_{1}d_{2}x^{i}\, d_{1}x^{j}d_{2}x^{k}+g_{ik,jl}\, d_{1}x^{i}d_{1}x^{j}\, d_{2}x^{k}d_{2}x^{l}\label{eq:ddgamma}\end{equation}
where $\Gamma_{ijk}=(g_{ij,k}+g_{ik,j}-g_{jk,i})/2$ are the Christoffel
symbols (components of the connection). If we look at (\ref{eq:ddgamma})
as a function of $d_{1}d_{2}x^{i}$'s and evaluate it at its critical
point \[
d_{1}d_{2}x^{i}=-\Gamma_{jk}^{i}\, d_{1}x^{i}d_{2}x^{j},\]
we get\begin{equation}
-\frac{1}{2}R_{ijkl}\, d_{1}x^{i}d_{1}x^{j}\, d_{2}x^{k}d_{2}x^{l}.\label{eq:R}\end{equation}
Geometrically we did the following: the space $(\Pi T)^{2}M=\HOM(\rodd{2},M)$
is fibred over the space $\Pi(T\oplus T)M$ of 1-jets of maps $\rodd{2}\map M$
(the additional coordinates in $(\Pi T)^{2}M$ are precisely $d_{1}d_{2}x^{i}$'s);
the critical points of $d_{1}d_{2}\gamma$ on the fibres of $(\Pi T)^{2}M\map\Pi(T\oplus T)M$
give us a section $\Pi(T\oplus T)M\map(\Pi T)^{2}M$; we use it to
pull back $d_{1}d_{2}\gamma$ and get a function on $\Pi(T\oplus T)M$,
i.e.~a section of $\bigwedge(T\oplus T)M$, (\ref{eq:R}).

Let us now look at the action of $\HOM(\rodd{2},\rodd{2})$ on the
worms $\gamma$, $d_{1}\gamma$, $d_{2}\gamma$ and $d_{1}d_{2}\gamma$.
The infinitesimal generators of the action were computed at the end
of Section \ref{sec:forms-worms}. We see directly that $R_{a}\gamma=0$
and $E_{b}^{a}\gamma=\delta_{b}^{a}\gamma$. As a result, $\gamma$,
$d_{1}\gamma$, $d_{2}\gamma$, $d_{1}d_{2}\gamma$ is the basis of
a representation of $\HOM(\rodd{2},\rodd{2})$, or better, they are
components of an equivariant map $\HOM(\rodd{2},M)\map V$, where
$V$ is certain irreducible 4-dimensional representation of $\HOM(\rodd{2},\rodd{2})$.
Geometrically, $V$ is the space of sections $\Gamma(S^{2}T^{*}\rodd{2})$.
Moreover, if $\tilde{\gamma}\in\Omega_{[2]}(M)$ is any worm satisfying
$R_{a}\tilde{\gamma}=0$ and $E_{b}^{a}\tilde{\gamma}=\delta_{b}^{a}\tilde{\gamma}$
then $\tilde{\gamma}=\tilde{g}_{ij}\, d_{1}x^{i}d_{2}x^{j}$ for some
tensor field $\tilde{g}\in\Gamma(S^{2}T^{*}M)$ (indeed: the condition
$E_{b}^{a}\tilde{\gamma}=\delta_{b}^{a}\tilde{\gamma}$ means that
$\tilde{\gamma}=\tilde{g}_{ij}\, d_{1}x^{i}d_{2}x^{j}+h_{i}\, d_{1}d_{2}x^{i}$
and $R_{a}\tilde{\gamma}=0$ gives $h_{i}=0$). In other words, an
equivariant map $\HOM(\rodd{2},M)\map V$ is the same as a section
of $S^{2}T^{*}M$.

Let us see what really happens. Let $g$ be a section of $S^{2}T^{*}M$.
Given any map $\phi:N\map M$ we get a section $\phi^{*}g$ of $S^{2}T^{*}N$.
If we take $N=\rodd{2}$, for any map $\rodd{2}\map M$ we get an
element of the vector space $V=\Gamma(S^{2}T^{*}\rodd{2})$, i.e.~we
have a map from $\HOM(\rodd{2},M)$ (the space of all maps $\rodd{2}\map M$)
to $V$. Our map $\HOM(\rodd{2},M)\map V$ is clearly equivariant.
Any tensor field $g\in S^{2}T^{*}M$ thus yields naturally an equivariant
map $\HOM(\rodd{2},M)\map V$, and as we noticed above, this correspondence
is a bijection. In other words, the sheaf of sections of $S^{2}T^{*}M$
is 2-representable (more information on approximations of covariant
tensor fields is in the last example in Section \ref{sec:Examples}).

\end{document}